# Superconvergence of the MINI mixed finite element discretization of the Stokes problem: An experimental study in 3D

Andrea Cioncolini[1] and Daniele Boffi[2,3]

[1] Department of Mechanical, Aerospace and Civil Engineering, University of Manchester, Oxford Road, Manchester M13 9PL, UK; andrea.cioncolini@manchester.ac.uk

[2] Computer, Electrical and Mathematical Science and Engineering Division, King Abdullah University of Science and Technology (KAUST), Thuwal 23955-6900, Kingdom of Saudi Arabia; daniele.boffi@kaust.edu.sa

[3] Dipartimento di Matematica ''F. Casorati'', University of Pavia, Via Ferrata 1, I-27100 Pavia, Italy

**Abstract**: Stokes flows are a type of fluid flow where convective forces are small in comparison with viscous forces, and momentum transport is entirely due to viscous diffusion. Besides being routinely used as benchmark test cases in numerical fluid dynamics, Stokes flows are relevant in several applications in science and engineering including porous media flow, biological flows, microfluidics, microrobotics, and hydrodynamic lubrication. The present study concerns the discretization of the equations of motion of Stokes flows in three dimensions utilizing the MINI mixed finite element, focusing on the superconvergence of the method which was investigated with numerical experiments using five purpose-made benchmark test cases with analytical solution. Despite the fact that the MINI element is only linearly convergent according to standard mixed finite element theory, a recent theoretical development proves that, for structured meshes in two dimensions, the pressure superconverges with order $O(h^{3/2})$, as well as the linear part of the computed velocity with respect to the piecewise-linear nodal interpolation of the exact velocity. The numerical experiments documented herein suggest a more general validity of the superconvergence in pressure, possibly to unstructured tetrahedral meshes and even up to quadratic convergence which was observed with one test problem, thereby indicating that there is scope to further extend the available theoretical results on convergence.

## 1. Introduction

In fluid dynamics, the expressions *Stokes flow* and, synonymously, *creeping flow* are employed when referring to fluid flows where inertia is small in comparison with viscous and pressure forces. The equation of motion for Stokes flows is the asymptotic limiting form of the Navier-Stokes equation for fluid dynamics when the Reynolds number becomes small, which in turn can be the result of a low flow velocity, a small characteristic flow length, a highly viscous fluid, or a combination thereof. In the limit of small Reynolds number values, momentum transport is only due to viscous diffusion whilst convection becomes inconsequential: the non-linear convection term in the full Navier-Stokes equation can therefore be neglected thereby linearizing the equation. Clearly, this simplification makes the equation of motion for Stokes flows more amenable to analytical treatment than the full Navier-Stokes equation, and indeed analytical solutions have been produced and documented for several Stokes flow problems. Consequently, Stokes flows have become essential and routinely used benchmark test cases in numerical fluid dynamics. On the practical side, Stokes flows are relevant in several applications in science and engineering including porous media flow [1], biological flows [2], microfluidics [3], microrobotics [4], and hydrodynamic lubrication [5].

In this study, we consider the discretization of the three-dimensional Stokes problem with the MINI mixed finite element, which we will refer to as *Stokes-MINI* problem, focusing in particular on experimentally investigating the superconvergence of the method. The theory of saddle point problems, developed by Babuška [6] and Brezzi [7], can be used as the foundation to inform the theoretical analysis of the Stokes problem discretization with mixed finite elements. The MINI mixed finite element, originally introduced by Arnold, Brezzi and Fortin [8] for the discretization of the Stokes problem, approximates the velocity space by continuous, piecewise-linear polynomials plus bubbles and the pressure space by continuous, piecewise-linear polynomials. The standard mixed finite element theory [9] shows that Stokes-MINI is stable and that quasi-optimal error estimates are satisfied. In particular, linear convergence $O(h)$ for both velocity and pressure can be easily proved. Notwithstanding the mismatch in the approximating properties of the finite element spaces employed (the velocity space is linearly convergent whereas the pressure space would allow for second order convergence), the Stokes-MINI discretization is fairly common because of its simplicity.

Despite the fact that the standard mixed finite element theory only guarantees $O(h)$ convergence for Stokes-MINI, Eichel, Tobiska and Xie [10] recently proved that the pressure and the linear part of the computed velocity with respect to the linear nodal interpolant of the exact velocity superconverge with order $O(h^{3/2})$. Presently, this superconvergence result is only proved in two dimensions on three-directional triangular meshes, which are structured triangular meshes obtained from a structured rectangular mesh when each rectangle is divided into two triangles using one of the rectangle diagonals. The existing superconvergence theory does not cover unstructured triangular meshes in two dimensions and unstructured tetrahedral meshes in three dimensions, thereby creating the motivation for the present work. Our main objective was to experimentally investigate $O(h^{3/2})$ superconvergence for the three-dimensional Stokes-MINI problem on a selection of five benchmark test cases with analytical solution, using unstructured tetrahedral meshes. To the best of our knowledge, no such experimental assessment for the three-dimensional Stokes problem is presently documented in the open literature. The test cases, which have been specifically designed for use here, represent enclosed flows defined in the unit cube, with Dirichlet-type boundary conditions for the velocity along the domain boundary. Two of these test cases, in particular, generalize the well-known lid-driven cavity flow problem to three dimensions.

Our experimental results suggest that the validity of the $O(h^{3/2})$ superconvergence in pressure goes beyond what covered by the existing theory, possibly to unstructured tetrahedral meshes and even up to $O(h^2)$ superconvergence, which we observed in one test problem. Conversely, we did not observe any superconvergence in velocity in our three-dimensional test cases, possibly indicating that $O(h^{3/2})$ superconvergence in velocity is restricted to the two-dimensional case. Our results also indicate that the polynomial bubble function, which is implemented to enrich the discrete velocity space so as to stabilize the MINI finite element discretization, generally has the effect of improving the quality of the velocity approximation but deteriorates the local mass conservation. The work documented herein represents a follow-on of our previous study [11], where we experimentally investigated $O(h^{3/2})$ superconvergence of Stokes-MINI on two-dimensional unstructured triangular meshes.

The rest of this paper is organized as follows: Section 2 provides the necessary theoretical background on the Stokes problem; the benchmark test cases with analytical solution are presented in Section 3; Section 4 describes the numerical methodology; the results of the numerical experiments are presented and discussed in Section 5; whilst the concluding remarks are summarized in Section 6.

## 2. Theoretical Background

### 2.1. The Stokes Problem

The strong formulation of the Stokes problem considered here reads as follows: find a three-dimensional velocity vector field $\underline{u}(x,y,z) = \big(u_x(x,y,z), u_y(x,y,z), u_z(x,y,z)\big)$, with $u_x, u_y, u_z \in C^2(\Omega) \cap C^0(\overline{\Omega})$, and a three-dimensional scalar pressure field $P(x,y,z) \in C^1(\Omega)$ such that:

$$-\mu \, \Delta \underline{u} + \underline{\nabla} P = \rho \underline{f} \quad in \ \Omega, \tag{1}$$

$$div\,(\underline{u}) = 0 \quad in \ \Omega, \tag{2}$$

$$\underline{u} = \underline{0} \quad on \ \partial\Omega, \tag{3}$$

where $\Omega \subset \mathbb{R}^3$ is a bounded and connected polyhedral domain in the space with boundary $\partial\Omega$, $\rho$ and $\mu$ are the (constant) fluid density and viscosity, $\Delta$, $\underline{\nabla}$ and $div$ are the Laplace, the gradient and the divergence operators, and $\underline{f} = (f_x, f_y, f_z)$ with $f_x, f_y, f_z \in C^0(\Omega)$ is the external force field. On the practical side, the Stokes problem in Eqs. (1) and (2) can effectively model the steady-state flow of incompressible and isothermal Newtonian fluids at low values of the Reynolds number. Even though we restrict our attention to Stokes flows with constant density and viscosity, the generalization to variable density and viscosity is straightforward. As usual, we consider homogeneous Dirichlet boundary condition for the velocity in Eq. (3), but it is understood that the present problem formulation can be extended to the non-homogeneous case: upon changing the variable and modifying the right-hand side of Eqs. (1) and (2), in fact, the non-homogeneous problem can be reformulated as a homogeneous one [9]. On the practical side, the Dirichlet boundary condition for the velocity is quite versatile and allows specifying inflows and outflows as well as adapting to the interaction of the fluid with solid boundaries. On the other hand, a Dirichlet boundary condition for the velocity imposed throughout the boundary only allows determining the pressure field up to an arbitrary additive constant [9], meaning that if a solution $(u_x, u_y, u_z, P)$ of the Stokes problem in Eqs. (1)-(3) exists then this is not unique, because $(u_x, u_y, u_z, P + c)$ will also be a solution for any constant $c \in \mathbb{R}$. Following common practice, we restored uniqueness by imposing a null mean value of the pressure field over the entire domain $\Omega$ as follows:

$$\frac{1}{|\Omega|} \int_\Omega P \, d\Omega = 0. \tag{4}$$

*2.2. Variational Formulation of the Stokes Problem*

The variational formulation of the Stokes problem is obtained in two steps [9]. First, Eqs. (1) and (2) are multiplied by test functions $\underline{v} = (v_x, v_y, v_z)$ and $q$, respectively, and then they are integrated over the domain $\Omega$. Then, after integration by parts and considering the boundary conditions, the variational formulation of the Stokes problem reads: find a velocity field $\underline{u}(x, y, z) = \big(u_x(x, y, z), u_y(x, y, z), u_z(x, y, z)\big)$, with $u_x, u_y, u_z \in H_0^1(\Omega)$ and a pressure function $P(x, y, z) \in L_0^2(\Omega)$ such that:

$$a(\underline{u}, \underline{v}) + b(\underline{v}, P) = F(\underline{v}) \quad \forall \, \underline{v} \in H_0^1(\Omega)^3, \tag{5}$$

$$b(\underline{u}, q) = 0 \quad \forall \, q \in L_0^2(\Omega), \tag{6}$$

with the following definitions for the bilinear and linear forms:

$$a(\underline{u}, \underline{v}) = \mu \int_\Omega \underline{\nabla}\,\underline{u} : \underline{\nabla}\,\underline{v} \, d\Omega = \mu \int_\Omega \underline{\nabla} u_x \cdot \underline{\nabla} v_x \, d\Omega + \mu \int_\Omega \underline{\nabla} u_y \cdot \underline{\nabla} v_y \, d\Omega + \mu \int_\Omega \underline{\nabla} u_z \cdot \underline{\nabla} v_z \, d\Omega, \tag{7}$$

$$b(\underline{v}, P) = -\int_\Omega div(\underline{v}) P \, d\Omega, \tag{8}$$

$$F(\underline{v}) = \rho \int_\Omega \underline{f} \cdot \underline{v} \, d\Omega, \tag{9}$$

where $\underline{f} = (f_x, f_y, f_z)$ with $f_x, f_y, f_z \in L^2(\Omega)$ is the external force field, $H_0^1(\Omega)$ is the Sobolev space of square-integrable functions vanishing on the domain boundary in the sense of traces and that have square-integrable first weak derivatives, whilst $L_0^2(\Omega)$ is the space of square-integrable functions with vanishing mean. From the standard theory of mixed problems [9], the Stokes problem in Eqs. (5) and (6) admits a unique solution $(\underline{u}, P) \in H_0^1(\Omega)^3 \times L_0^2(\Omega)$.

*2.3. Finite Element Approximation of the Stokes Problem*

The Galerkin approximation of the Stokes problem in Eqs. (5) and (6) is obtained when the problem is formulated by considering finite-dimensional linear subspaces $X_h(\Omega) \subset H_0^1(\Omega)$ and $M_h(\Omega) \subset L_0^2(\Omega)$, depending on a discretization parameter $h$, and which approximate the Hilbert spaces $H_0^1(\Omega)$ and $L_0^2(\Omega)$, so that we are led to the discrete problem: find $(\underline{u}_h, P_h) \in X_h^3 \times M_h$ such that:

$$a(\underline{u}_h, \underline{v}_h) + b(\underline{v}_h, P_h) = F(\underline{v}_h) \quad \forall \, \underline{v}_h \in X_h(\Omega)^3, \tag{10}$$

$$b(\underline{u}_h, q_h) = 0 \quad \forall \, q_h \in M_h(\Omega). \tag{11}$$

The discrete Stokes problem in Eqs. (10) and (11) is uniquely solvable if the discrete linear subspaces $X_h$ and $M_h$ satisfy the following discrete inf-sup condition [9]:

$$\exists \, \beta > 0: \quad \inf_{0 \neq q_h \in M_h} \sup_{0 \neq \underline{v}_h \in X_h^3} \frac{b(\underline{v}_h, q_h)}{\| \underline{v}_h \|_X \, \| q_h \|_M} \geq \beta, \tag{12}$$

where $\beta$ is a constant independent of $h$. On the practical side, the discrete inf-sup condition in Eq. (12) is a compatibility condition between the spaces $X_h$ and $M_h$ to inform the development of uniquely solvable numerical schemes for the discrete Stokes problem.

*2.4. The MINI Element for the Stokes Problem*

When the Stokes problem is approximated using mixed finite elements, the discrete linear subspaces $X_h$ and $M_h$ comprise piecewise polynomial functions defined on a discretization $\Omega_h$ of the domain $\Omega$. The domain discretizations of interest here are conformal tetrahedral meshes $\mathcal{T}_h$ where any two tetrahedra share at most one vertex, one edge, or one face. Since the domain $\Omega$ is bounded, connected and polyhedral, its closure $\overline{\Omega}$ can be divided into tetrahedra $T$ that form a mesh $\mathcal{T}_h$ that wholly covers $\Omega$, so that $\Omega_h = \mathcal{T}_h \equiv \Omega$. The discretization parameter $h$, in particular, corresponds the maximum diameter of the tetrahedra in the mesh. Several mixed finite elements have been proposed for numerical fluid dynamics applications and the choice of the "best element" depends on various aspects, including the quantities of interest in the problem under consideration and the coding infrastructure. Finite elements based on discontinuous pressures generally provide a more localized mass conservation, since in Eq. (11) it is possible to take test functions $q_h$ supported in a single element. This comes at the price of a larger discrete pressure space; in this case direct solvers can be efficiently used in combination with penalty methods. On the other hand, finite elements based on continuous pressures might be preferred by programmers more used to standard finite elements; iterative solvers are more efficient in this case since the size of the matrix and the computational cost are reduced. The method we are considering belongs to the latter category. Basically, distinct mixed finite elements differ in the local order and in the global regularity of the polynomials employed, and are identified with the symbol $[\mathbb{P}_k]^d / \mathbb{P}_m$ where $d$ is the number of space dimensions (typically 2 or 3) whilst $k$ and $m$ and the orders of the polynomials adopted for the velocity and pressure spaces.

The MINI element was originally introduced by Arnold, Brezzi and Fortin [8] and is based on continuous-piecewise-linear polynomials enriched with local polynomial bubble functions for the discrete velocity space $X_h$ and continuous-piecewise-linear polynomials for the discrete pressure space $M_h$, which corresponds to $[\mathbb{P}_{1+b}]^2/\mathbb{P}_1$ in two dimensions and to $[\mathbb{P}_{1+b}]^3/\mathbb{P}_1$ in three dimensions. The starting point of the authors was the $[\mathbb{P}_1]^2/\mathbb{P}_1$ pair for two-dimensional problems and its extension $[\mathbb{P}_1]^3/\mathbb{P}_1$ for three-dimensional problems, which would be computationally quite convenient but that are unfortunately unstable because they do not satisfy the discrete inf-sup condition in Eq. (12) [9]. The MINI element was therefore conceived as a stabilized version of the $[\mathbb{P}_1]^2/\mathbb{P}_1$ pair and of its three-dimensional extension $[\mathbb{P}_1]^3/\mathbb{P}_1$ where local functions, named *bubble functions*, are included to enrich the discrete velocity space thereby stabilizing the scheme. Further details on the two-dimensional MINI element can be found in [8,9,11]. In the three-dimensional case of interest here, the bubble function is a fourth-degree polynomial locally defined in each tetrahedron as the product of the linear nodal basis functions (so-called barycentric coordinates) $\varphi_1, \varphi_2, \varphi_3, \varphi_4$ of the tetrahedron itself (see Eq. (14) where $v_{h|T}$ refers to one of the three components of the discrete velocity):

$$X_h = \{v_h \in C^0(\overline{\Omega}): v_{h|T} \in \mathbb{P}_{1+b} \ \forall T \in \mathcal{T}_h\}, \tag{13}$$

$$v_{h|T} = a_h + b_h x + c_h y + d_h z + e_h \ \varphi_1(x,y,z)\varphi_2(x,y,z)\varphi_3(x,y,z)\varphi_4(x,y,z), \tag{14}$$

$$M_h = \{q_h \in C^0(\overline{\Omega}): q_{h|T} \in \mathbb{P}_1 \ \forall T \in \mathcal{T}_h\}, \tag{15}$$

$$q_{h|T} = a_h + b_h x + c_h y + d_h z, \tag{16}$$

where $a_h, b_h, c_h, d_h, e_h$ are constants. As can be noted from inspecting Eq. (14), with the MINI element the discrete velocity components are expressed as the sum of a linear polynomial plus the bubble function. The computed velocity $\underline{u}_h$ can therefore be expressed as the sum of a piecewise-linear part $\underline{u}_{hl}$ plus a bubble part $\underline{u}_{hb}$ as follows:

$$\underline{u}_h = \underline{u}_{hl} + \underline{u}_{hb}, \tag{17}$$

As noted by Verfürth [12], Bank and Welfert [13,14], Kim and Lee [15] and Russo [16], in practical applications the piecewise-linear velocity $\underline{u}_{hl}$ is often observed to approximate the exact velocity $\underline{u}$ better than the complete computed velocity $\underline{u}_h$ itself. This suggests that the bubble part $\underline{u}_{hb}$, which is crucial to stabilize the formulation, does not significantly contribute to reduce the error in the velocity. Accordingly, a-posteriori error estimators are often constructed using the piecewise-linear velocity $\underline{u}_{hl}$ of the Stokes-MINI discretization [12].

The standard mixed finite element theory [9] guarantees the stability and the linear convergence of Stokes-MINI:

$$\| \underline{u} - \underline{u}_h \|_{H^1} + \| P - P_h \|_{L^2} \leq C\, h \left( \| \underline{u} \|_{H^2} + \| P \|_{H^1} \right), \tag{18}$$

where $C$ is a positive constant independent of $h$, $\|*\|_{H^1}$ and $\|*\|_{L^2}$ are the usual norms in $H^1(\Omega)$ and $L^2(\Omega)$, and provided the exact solution $(\underline{u}, P) \in H^2(\Omega)^2 \times H^1(\Omega)$ in two dimensions or $(\underline{u}, P) \in H^2(\Omega)^3 \times H^1(\Omega)$ in three dimensions.

Despite the fact that the standard mixed finite element theory only guarantees $O(h)$ convergence for the Stokes-MINI problem, Eichel, Tobiska and Xie [10] recently proved the following $O(h^{3/2})$ superconvergence result in two dimensions on three-directional triangular meshes, that is, structured triangular meshes obtained from a rectangular mesh after decomposition of each rectangle into two triangles (using one of the rectangle diagonals):

**Theorem 1.** *With reference to the two-dimensional Stokes-MINI problem, assume the following regularity of the exact solution $(\underline{u}, P) \in H^3(\Omega)^2 \times H^2(\Omega)$ and consider a three-directional triangulation $\mathcal{T}_h$, then:*

$$\left( |\underline{u}_{hl} - i_h \underline{u}|_{H^1}^2 + \| P_h - j_h P \|_{L^2}^2 \right)^{1/2} \leq C\, h^{3/2} \left( \| \underline{u} \|_{H^3} + \| P \|_{H^2} \right), \tag{19}$$

*where $C$ is a positive constant independent of $h$, $|*|_{H^1}$ is the usual semi-norm in $H^1(\Omega)$, and $(i_h \underline{u}, j_h P)$ denotes the vertex-based piecewise-linear nodal interpolation of the exact solution $(\underline{u}, P)$.*

The proof can be found in [10]. Theorem 1 shows a faster rate of convergence $O(h^{3/2})$ when the piecewise-linear computed velocity $\underline{u}_{hl}$ is compared to the piecewise-linear nodal interpolation of the exact velocity $i_h \underline{u}$, and the discrete pressure $P_h$ is compared to the piecewise-linear nodal interpolation of the exact pressure $j_h P$. Under the same assumptions, the superconvergence in pressure can be further generalized as follows [10]:

$$\| P - P_h \|_{L^2} \leq \| P_h - j_h P \|_{L^2} + \| j_h P - P \|_{L^2} \leq C \left( h^{3/2} + h^2 \right) \left( \| \underline{u} \|_{H^3} + \| P \|_{H^2} \right), \tag{20}$$

which guarantees $O(h^{3/2})$ superconvergence of the computed pressure $P_h$ to the exact pressure $P$. Presently, this superconvergence results is only proved on three-directional triangular meshes in two-dimensions; unstructured meshes and the three-dimensional case are not covered.

## 3. Benchmark Test Cases with Analytical Solution

The five test cases with analytical solution specifically designed for use here and described below represent enclosed flows defined in the unit cube with the velocity specified throughout the boundary and vanishing mean of the pressure over the domain. In particular, in three test cases (#1, #2, and #3) the boundary condition for the velocity is homogeneous, whilst two test cases (#4 and #5) have non-homogeneous boundary conditions. For these latter, only the velocity component tangential to the boundary is non-zero, whereas the normal component is zero. This way, it is possible to deal with a more general non-homogeneous boundary condition without having the burden of numerically handling inflows and outflows. As it is well known, incompressible flow problems with inflow/outflow are particularly challenging to solve because of the difficulty of enforcing global mass conservation at the discrete level. This is the reason why, when feasible, it is preferred to reformulate inflow/outflow boundary conditions in incompressible flow problems [17]. The test cases #4 and #5, therefore, can be regarded as three-dimensional generalizations of the lid-driven cavity flow problem, which is one of the most widely used benchmark validation cases in two-dimensional numerical fluid dynamics.

*3.1. Test Problem #1*

Consider the Stokes problem in Eqs. (1) and (2) defined on the unit cube domain $\Omega = (0,1) \times (0,1) \times (0,1)$ with external force given as follows:

$$\begin{aligned}\rho f_x = -\mu[&(2 - 12x + 12x^2)(2y - 6y^2 + 4y^3)(2z - 6z^2 + 4z^3) \\ &+ (x^2 - 2x^3 + x^4)(-12 + 24y)(2z - 6z^2 + 4z^3) \\ &+ (x^2 - 2x^3 + x^4)(2y - 6y^2 + 4y^3)(-12 + 24z)] + 0.01,\end{aligned} \quad (21)$$

$$\begin{aligned}\rho f_y = -\mu[&(-12 + 24x)(y^2 - 2y^3 + y^4)(2z - 6z^2 + 4z^3) \\ &+ (2x - 6x^2 + 4x^3)(2 - 12y + 12y^2)(2z - 6z^2 + 4z^3) \\ &+ (2x - 6x^2 + 4x^3)(y^2 - 2y^3 + y^4)(-12 + 24z)] + 0.01,\end{aligned} \quad (22)$$

$$\begin{aligned}\rho f_z = 2\mu[&(-12 + 24x)(2y - 6y^2 + 4y^3)(z^2 - 2z^3 + z^4) \\ &+ (2x - 6x^2 + 4x^3)(-12 + 24y)(z^2 - 2z^3 + z^4) \\ &+ (2x - 6x^2 + 4x^3)(2y - 6y^2 + 4y^3)(2 - 12z + 12z^2)] - 0.01,\end{aligned} \quad (23)$$

and with null velocity on all domain boundary:

$$u_x = u_y = u_z = 0 \ \ on \ \partial\Omega; \quad (24)$$

the corresponding exact solution is:

$$u_x = (x^2 - 2x^3 + x^4)(2y - 6y^2 + 4y^3)(2z - 6z^2 + 4z^3), \quad (25)$$

$$u_y = (2x - 6x^2 + 4x^3)(y^2 - 2y^3 + y^4)(2z - 6z^2 + 4z^3), \quad (26)$$

$$u_z = -2(2x - 6x^2 + 4x^3)(2y - 6y^2 + 4y^3)(z^2 - 2z^3 + z^4), \quad (27)$$

$$P = 0.01(x + y + z - 1.5). \quad (28)$$

*3.2. Test Problem #2*

Consider the Stokes problem in Eqs. (1) and (2) defined on the unit cube domain $\Omega = (0,1) \times (0,1) \times (0,1)$ with external force given as follows:

$$\rho f_x = -4\pi^2 \mu [3 \cos(2\pi x) - 2] \sin(2\pi y) \sin(2\pi z) - 2\pi \sin(2\pi x), \tag{29}$$

$$\rho f_y = -4\pi^2 \mu \sin(2\pi x)[3 \cos(2\pi y) - 2] \sin(2\pi z) - 2\pi \sin(2\pi y), \tag{30}$$

$$\rho f_z = 8\pi^2 \mu \sin(2\pi x) \sin(2\pi y)[3 \cos(2\pi z) - 2] - 2\pi \sin(2\pi z), \tag{31}$$

and with null velocity on all domain boundary:

$$u_x = u_y = u_z = 0 \; on \; \partial\Omega; \tag{32}$$

the corresponding exact solution is:

$$u_x = [1 - \cos(2\pi x)] \sin(2\pi y) \sin(2\pi z), \tag{33}$$

$$u_y = \sin(2\pi x)[1 - \cos(2\pi y)] \sin(2\pi z), \tag{34}$$

$$u_z = -2 \sin(2\pi x) \sin(2\pi y)[1 - \cos(2\pi z)], \tag{35}$$

$$P = \cos(2\pi x) + \cos(2\pi y) + \cos(2\pi z). \tag{36}$$

*3.3. Test Problem #3*

Consider the Stokes problem in Eqs. (1) and (2) defined on the unit cube domain $\Omega = (0,1) \times (0,1) \times (0,1)$ with external force field given as follows:

$$\begin{aligned}\rho f_x = -\mu e^x [&(2 - 8x + x^2 + 6x^3 + x^4)(2y - 6y^2 + 4y^3)(2z - 6z^2 + 4z^3) \\ &+ (x^2 - 2x^3 + x^4)(-12 + 24y)(2z - 6z^2 + 4z^3) \\ &+ (x^2 - 2x^3 + x^4)(2y - 6y^2 + 4y^3)(-12 + 24z)] + 0.01yz,\end{aligned} \tag{37}$$

$$\begin{aligned}\rho f_y = \mu e^x [&(-8 + 2x + 18x^2 + 4x^3)(y^2 - 2y^3 + y^4)(2z - 6z^2 + 4z^3) \\ &+ (2x - 6x^2 + 4x^3)(2 - 12y + 12y^2)(2z - 6z^2 + 4z^3) \\ &+ (2x - 6x^2 + 4x^3)(y^2 - 2y^3 + y^4)(-12 + 24z)] + 0.01xz,\end{aligned} \tag{38}$$

$$\begin{aligned}\rho f_z = \mu e^x [&(2 - 8x + x^2 + 6x^3 + x^4)(2y - 6y^2 + 4y^3)(z^2 - 2z^3 + z^4) \\ &+ (x^2 - 2x^3 + x^4)(-12 + 24y)(z^2 - 2z^3 + z^4) \\ &+ (x^2 - 2x^3 + x^4)(2y - 6y^2 + 4y^3)(2 - 12z + 12z^2)] + 0.01xy,\end{aligned} \tag{39}$$

and with null velocity on all domain boundary:

$$u_x = u_y = u_z = 0 \; on \; \partial\Omega; \tag{40}$$

the corresponding exact solution is:

$$u_x = e^x(x^2 - 2x^3 + x^4)(2y - 6y^2 + 4y^3)(2z - 6z^2 + 4z^3), \tag{41}$$

$$u_y = -e^x(2x - 6x^2 + 4x^3)(y^2 - 2y^3 + y^4)(2z - 6z^2 + 4z^3), \tag{42}$$

$$u_z = -e^x(x^2 - 2x^3 + x^4)(2y - 6y^2 + 4y^3)(z^2 - 2z^3 + z^4), \tag{43}$$

$$P = 0.01(xyz - 0.125). \tag{44}$$

*3.4. Test Problem #4*

Consider the Stokes problem in Eqs. (1) and (2) defined on the unit cube domain $\Omega = (0,1) \times (0,1) \times (0,1)$ with external force field given as follows:

$$\begin{aligned}\rho f_x = &-\mu[(2 - 12x + 12x^2)(2y - 6y^2 + 4y^3)(-z + 2z^3) \\ &+ (x^2 - 2x^3 + x^4)(-12 + 24y)(-z + 2z^3) \\ &+ (x^2 - 2x^3 + x^4)(2y - 6y^2 + 4y^3)12z] \\ &+ (1 - 6x + 6x^2)(y - 3y^2 + 2y^3)(z - 3z^2 + 2z^3),\end{aligned} \tag{45}$$

$$\begin{aligned}\rho f_y = &-\mu[(-12 + 24x)(y^2 - 2y^3 + y^4)(-z + 2z^3) \\ &+ (2x - 6x^2 + 4x^3)(2 - 12y + 12y^2)(-z + 2z^3) \\ &+ (2x - 6x^2 + 4x^3)(y^2 - 2y^3 + y^4)12z] \\ &+ (x - 3x^2 + 2x^3)(1 - 6y + 6y^2)(z - 3z^2 + 2z^3),\end{aligned} \tag{46}$$

$$\begin{aligned}\rho f_z = &\mu[(-12 + 24x)(2y - 6y^2 + 4y^3)(-z^2 + z^4) \\ &+ (2x - 6x^2 + 4x^3)(-12 + 24y)(-z^2 + z^4) \\ &+ (2x - 6x^2 + 4x^3)(2y - 6y^2 + 4y^3)(-2 + 12z^2)] \\ &+ (x - 3x^2 + 2x^3)(y - 3y^2 + 2y^3)(1 - 6z + 6z^2),\end{aligned} \tag{47}$$

and with null velocity on all domain boundary except along the top face of the domain (note that only the tangential velocity component is non-zero whereas the normal velocity component is zero, so that there is no inflow nor outflow at the top face of the domain boundary):

$$u_x(x, y, 0) = u_y(x, y, 0) = u_z(x, y, 0) = 0, \tag{48}$$

$$u_x(x, y, 1) = (x^2 - 2x^3 + x^4)(2y - 6y^2 + 4y^3), \tag{49}$$

$$u_y(x, y, 1) = (2x - 6x^2 + 4x^3)(y^2 - 2y^3 + y^4),$$

$$u_z(x, y, 1) = 0,$$

$$u_x(x, 0, z) = u_y(x, 0, z) = u_z(x, 0, z) = 0, \tag{50}$$

$$u_x(x, 1, z) = u_y(x, 1, z) = u_z(x, 1, z) = 0, \quad (51)$$

$$u_x(0, y, z) = u_y(0, y, z) = u_z(0, y, z) = 0, \quad (52)$$

$$u_x(1, y, z) = u_y(1, y, z) = u_z(1, y, z) = 0; \quad (53)$$

the corresponding exact solution is:

$$u_x = (x^2 - 2x^3 + x^4)(2y - 6y^2 + 4y^3)(-z + 2z^3), \quad (54)$$

$$u_y = (2x - 6x^2 + 4x^3)(y^2 - 2y^3 + y^4)(-z + 2z^3), \quad (55)$$

$$u_z = -(2x - 6x^2 + 4x^3)(2y - 6y^2 + 4y^3)(-z^2 + z^4), \quad (56)$$

$$P = (x - 3x^2 + 2x^3)(y - 3y^2 + 2y^3)(z - 3z^2 + 2z^3). \quad (57)$$

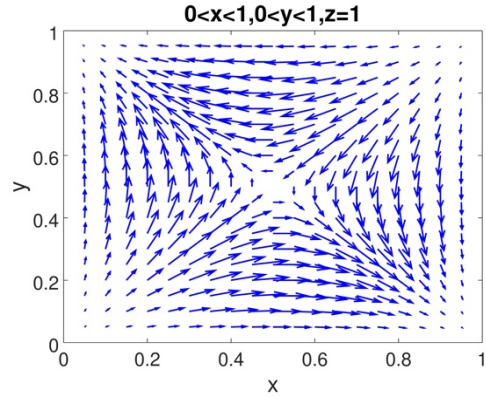

**Figure 1.** Tangential velocity vector field along the top face of the domain for Test Problem #4.

With the present test problem, the velocity field is driven by the body force specified in Eqs. (45)–(47) and by the non-uniform shear acting on the top face of the domain boundary in Eq. (49), so that this flow problem can be regarded as a three-dimensional generalization of the lid-driven cavity flow. Notably, there are no velocity singularities at the corners or edges of the top face of the domain. The tangential velocity vector field along the top face of the domain is provided in Figure 1.

*3.5. Test Problem #5*

Consider the Stokes problem in Eqs. (1) and (2) defined on the unit cube domain $\Omega = (0,1) \times (0,1) \times (0,1)$ with external force field given as follows:

$$\rho f_x = (3\mu\pi^2 - \pi) \sin(\pi x) \cos(\pi y) \cos(\pi z), \quad (58)$$

$$\rho f_y = (3\mu\pi^2 - \pi) \cos(\pi x) \sin(\pi y) \cos(\pi z), \quad (59)$$

$$\rho f_z = -(6\mu\pi^2 + \pi) \cos(\pi x) \cos(\pi y) \sin(\pi z), \tag{60}$$

and with non-homogeneous Dirichlet boundary condition for the velocity on all domain boundaries (note that only the tangential velocity component is non-zero whereas the normal velocity component is zero, so that there is no inflow nor outflow at the domain boundary):

$$u_x(x, y, 0) = sin(\pi x) \, cos(\pi y),$$
$$u_y(x, y, 0) = cos(\pi x) \, sin(\pi y), \tag{61}$$
$$u_z(x, y, 0) = 0,$$

$$u_x(x, y, 1) = -sin(\pi x) \, cos(\pi y),$$
$$u_y(x, y, 1) = -cos(\pi x) \, sin(\pi y), \tag{62}$$
$$u_z(x, y, 1) = 0,$$

$$u_x(x, 0, z) = sin(\pi x) \, cos(\pi z)$$
$$u_y(x, 0, z) = 0, \tag{63}$$
$$u_z(x, 0, z) = -2 \, cos(\pi x) \, sin(\pi z),$$

$$u_x(x, 1, z) = -sin(\pi x) \, cos(\pi z)$$
$$u_y(x, 1, z) = 0, \tag{64}$$
$$u_z(x, 1, z) = 2 \, cos(\pi x) \, sin(\pi z),$$

$$u_x(0, y, z) = 0,$$
$$u_y(0, y, z) = sin(\pi y) \, cos(\pi z) \tag{65}$$
$$u_z(0, y, z) = -2 \, cos(\pi y) \, sin(\pi z),$$

$$u_x(1, y, z) = 0,$$
$$u_y(1, y, z) = -sin(\pi y) \, cos(\pi z), \tag{66}$$
$$u_z(1, y, z) = 2 \, cos(\pi y) \, sin(\pi z);$$

the corresponding exact solution is:

$$u_x = sin(\pi x) \, cos(\pi y) \, cos(\pi z), \tag{67}$$

$$u_y = cos(\pi x) \, sin(\pi y) \, cos(\pi z), \tag{68}$$

$$u_z = -2\cos(\pi x)\cos(\pi y)\sin(\pi z), \tag{69}$$

$$P = \cos(\pi x)\cos(\pi y)\cos(\pi z). \tag{70}$$

With the present test problem, the velocity field is driven by the body force specified in Eqs. (58)–(60) and by the non-uniform shear acting along the faces of the domain boundary in Eqs. (61)–(70), so that this flow problem can be regarded as a three-dimensional generalization of the lid-driven cavity flow. Tangential velocity vector fields along the faces of the domain are provided in Figure 2.

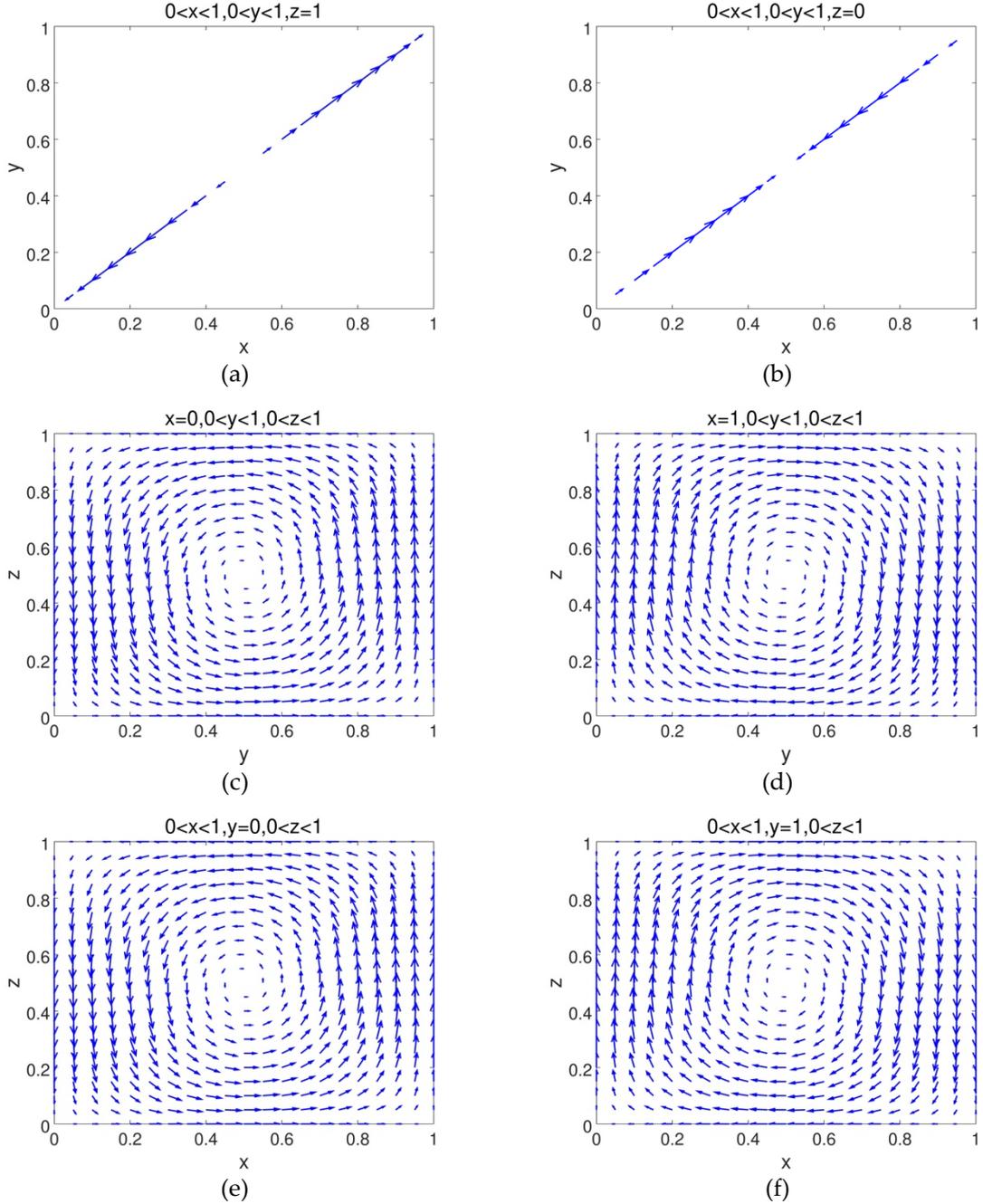

**Figure 2**. Tangential velocity vector fields along the faces of the domain for Test Problem #5: (a) top face, (b) bottom face, (c–f) side faces.

## 4. Numerical Methodology

We performed all calculations with the free software GNU Octave version 4.2.2 [18] (running under Linux Ubuntu 18.04.5 LTS), setting to 1 kg/ms the kinematic viscosity $\mu$ in Eq. (1). The mesh generation, the numerical integration, and the solution of the linear system are explained below with a level of detail sufficient to allow reproducing the presented results.

*4.1. Mesh generation*

The computational meshes were generated using the free software DistMesh (downloaded at: https://github.com/ionhandshaker/distmesh), a simple MATLAB and GNU Octave code for unstructured triangular and tetrahedral meshes [19]. DistMesh handles the computational domain using a signed distance function and generates the mesh by exploiting a physical analogy between a truss structure and a simplex mesh. In particular, DistMesh implements a suitable force-displacement relation for the beams and solves for the equilibrium of the truss structure: the forces move the nodes whilst the mesh topology is adjusted with the Delaunay triangulation algorithm. This results in a simple yet effective mesh generator that typically produces meshes of good quality [19].

Since all test problems considered here are defined on the same unit cube domain, they were solved using the same set of twelve meshes, whose main properties are summarized in Table 1.

**Table 1.** Main parameters of the tetrahedral meshes used in the calculations

| Mesh # | $h$ | No. Verts | No. Tets | Shape Ratio min/mean/max | Dihedral Angle (°) min/mean/max |
|---|---|---|---|---|---|
| 1 | 0.408 | 125 | 399 | 0.173/0.824/0.999 | 8.1/69.6/165.1 |
| 2 | 0.327 | 222 | 780 | 0.161/0.851/0.996 | 8.3/69.8/167.3 |
| 3 | 0.280 | 341 | 1313 | 0.067/0.871/0.998 | 3.5/69.8/174.9 |
| 4 | 0.251 | 505 | 2072 | 0.065/0.878/0.997 | 3.4/69.8/175.0 |
| 5 | 0.205 | 721 | 3090 | 0.030/0.883/0.998 | 1.5/69.9/177.7 |
| 6 | 0.169 | 1339 | 6106 | 0.036/0.896/0.999 | 1.9/69.9/177.2 |
| 7 | 0.142 | 2203 | 10534 | 0.048/0.898/0.999 | 2.4/69.9/176.2 |
| 8 | 0.120 | 4076 | 20333 | 0.011/0.902/0.999 | 0.6/69.9/179.1 |
| 9 | 0.100 | 9221 | 48100 | 0.007/0.904/1.000 | 0.3/69.9/179.4 |
| 10 | 0.084 | 12014 | 63817 | 0.008/0.908/0.999 | 0.4/70.0/179.4 |
| 11 | 0.061 | 32699 | 179325 | 0.006/0.909/1.000 | 0.2/70.0/179.6 |
| 12 | 0.048 | 50546 | 279408 | 0.004/0.910/1.000 | 0.2/70.0/179.7 |

All meshes were generated with uniform spacing. The number of vertices and tetrahedra in the meshes varied from 125 to 50,546 and from 399 to 279,408, respectively, whilst the mesh spacing parameter $h$ (corresponding to the length of the longest edge of all tetrahedra in the mesh) ranged from 0.408 to 0.048. The quality of the tetrahedra was assessed using the shape ratio $\rho$ (also called 'radius ratio' or 'aspect ratio'), which is defined as:

$$\rho = \frac{3\,r}{R}, \tag{71}$$

where $r$ and $R$ are the inradius (the radius of the smallest sphere that is tangent to all four faces of the tetrahedron) and the circumradius (the radius of the unique sphere that passes through all four vertices of the tetrahedron), respectively, and can be calculated as follows:

$$r = \frac{3\, V_{tet}}{A_1 + A_2 + A_3 + A_4}, \quad (72)$$

$$R = \frac{\sqrt{(a+b+c)(a+b-c)(a-b+c)(-a+b+c)}}{24\, V_{tet}}, \quad (73)$$

where $V_{tet}$ is the volume of the tetrahedron, $A_1, A_2, A_3, A_4$ are the areas of its faces, and $a, b, c$ are the products of the lengths of its opposite edges. The shape ratio, which is a relatively simple to compute yet effective and widely used shape measure for tetrahedral meshes [20,21], is dimensionless and bounded between 0 and 1: for a regular tetrahedron in which all four faces are equilateral triangles the shape ratio is 1, whilst it is 0 for a degenerate tetrahedron with zero volume.

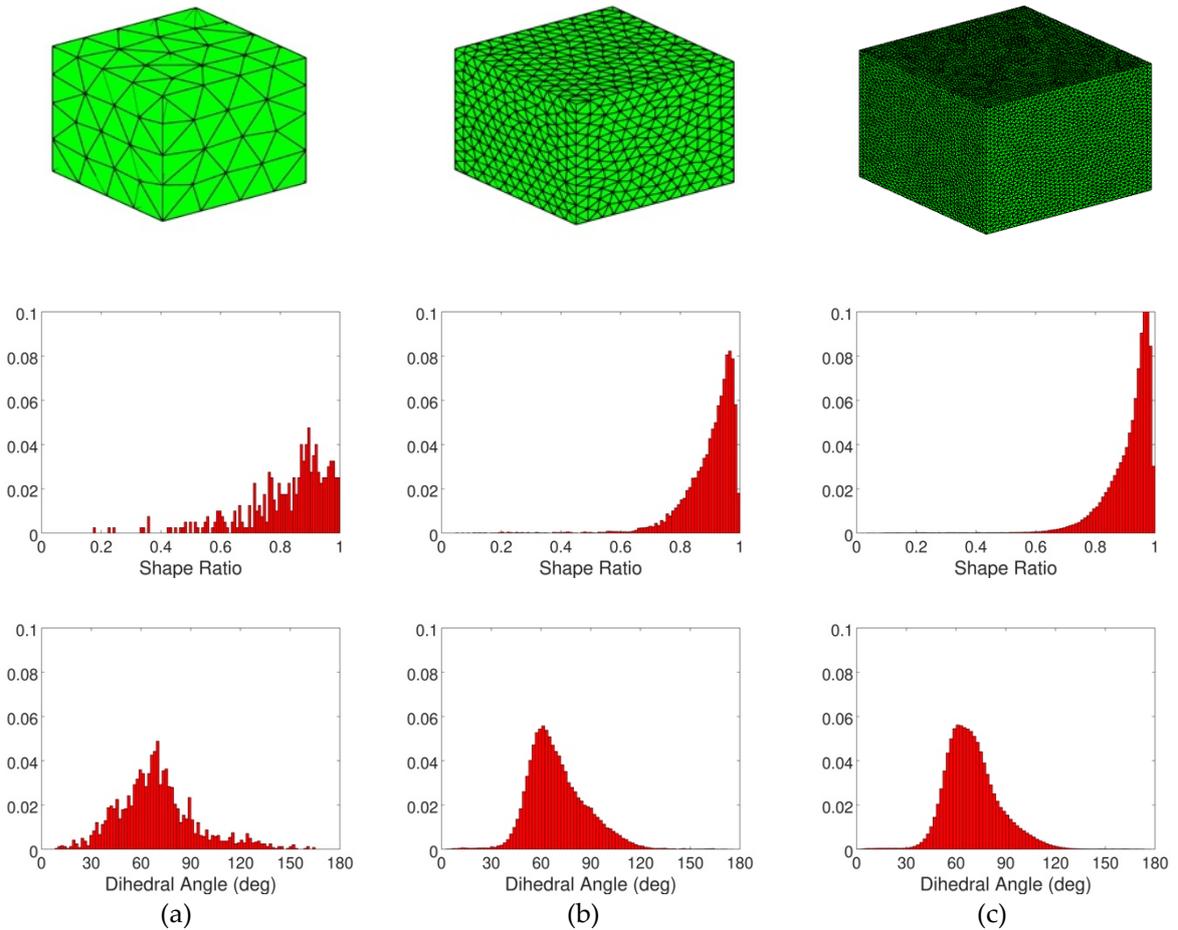

**Figure 3.** Mesh rendering and normalized histograms of shape ratios and dihedral angles for three representative meshes from Table 1: (a) mesh #1 (coarsest), (b) mesh #7 (intermediate), (c) mesh #12 (finest).

Following common practice, the quality of the tetrahedra was further assessed using the dihedral angles, which are the angles between adjoining triangular faces (a tetrahedron has six dihedral angles, one for each edge). Operatively, the dihedral angles were computed by taking the scalar products of the unit normal vectors to the triangular faces of the tetrahedron. The unit normal vectors, in turn, can be easily computed from the coordinates of the vertices. For a regular tetrahedron in which all four faces are equilateral triangles the dihedral angles are equal to ≈ 70.56°, whilst in a degenerate tetrahedron with zero volume at least one dihedral angle equals 0° or 180°.

The minimum, the maximum and the mean value of the shape ratios and dihedral angles of all meshes are provided in Table 1, whilst mesh renderings and normalized histograms are provided in Figure 3 for three representative meshes from Table 1 (the coarsest mesh #1, the intermediate mesh #7, and the finest mesh #12). As can be noted in Table 1, the mean values of the shape ratio vary in the range of 0.824–0.910 and those of the dihedral angle vary in the range of 69.6°–70.0°, which are indicative of good quality meshes. However, all meshes include a relatively small number of low-quality tetrahedra characterized by small values of the shape ratio and very small or very large dihedral angles. In particular, the quality of the worst tetrahedron in the mesh seems to degrade as the mesh is gradually refined. The low-quality tetrahedra in the present meshes are very few in number: referring to mesh #12 as a representative example, out of 279,408 total tetrahedra 2367 have shape ratio below 0.5 (0.85% of the total) and 160 have shape ratio below 0.1 (0.057% of the total). Further inspection indicated that these low-quality tetrahedra are mostly slivers: tetrahedra which have reasonable edge lengths and reasonable face areas but small volume because their vertices lie close to a plane.

Unfortunately, slivers are quite common in three-dimensional Delaunay triangulations, so that the present results are not surprising. As it is well known, slivers in finite element computations can negatively affect the discretization error and the conditioning of the stiffness matrix [22]. Even though post-processing techniques have been developed to improve the quality of tetrahedral meshes containing slivers [23,24], the meshes in Table 1 were considered appropriate for the scope of the present study, and were therefore not improved or otherwise modified. As stated previously, our main objective was to experimentally assess superconvergence. In this respect, it is the rate of change of the discretization error as the mesh is gradually refined what really matters, whereas the absolute value of the discretization error is not critical, and a slight worsening due to slivers would not be necessarily detrimental. For what concerns the conditioning of the stiffness matrix, solving the linear system posed no concern (see discussion later on), indicating that any worsening due to slivers did not go beyond what the preconditioner could handle. Moreover, it is not clear yet how numerous slivers need to be in a tetrahedral mesh in order to cause serious damage: as noted by Persson and Strang [19] specifically for DistMesh, a few isolated slivers might not be seriously detrimental after all.

*4.2. Numerical Integration*

In order to evaluate the components of the stiffness matrix, the components of the coefficient vector and the discretization errors, various functions need to be numerically integrated over the mesh tetrahedra. The numerical integration was carried out employing the quadrature formulas for tetrahedra developed by Witherden and Vincent [25] specifically for finite element applications. These quadrature formulas are symmetric, have degree of up to 10, have strictly positive weights, and all quadrature nodes are located inside the domain. For polynomial integrands of degree up to 10 the appropriate formula with the same degree as the integrand function was used, whilst polynomial integrands of degree higher than 10 and non-polynomial integrands were evaluated with the formula of highest degree (degree 10).

*4.3. Solution of the Linear System*

The linear system was numerically solved with the Generalized Minimum Residual Method (GMRES) preconditioned with incomplete LU factorization (ILU), using the GNU Octave built-in functions *gmres* and *ilu*. When dealing with symmetric and indefinite linear systems such as those considered here, Krylov subspace iterative methods such as GMRES are particularly effective [26], provided that an effective preconditioning strategy is implemented [27]. GMRES with ILU preconditioning, in particular, is an effective go-to technique that has performed well in a wide range of practical problems [28], thereby motivating its use here. Another popular method for the Stokes problem is the Minimum Residual Method (MINRES), which is potentially more efficient than GMRES because it takes advantage of the symmetry of the system matrix. GMRES was however preferred to MINRES because designing an efficient preconditioner for this latter requires more effort. Operatively, for the present calculations the drop-tolerance in the ILU factorization was set at $10^{-4}$, the target relative residual in GMRES was set equal to $10^{-9}$, and the initial guess in GMRES was the zero vector (the default in GNU Octave). Overall, the performance of the method was quite satisfactory: the number of GMRES iterations required to reach convergence for the five test cases ranged between 9–10 for the coarsest mesh (mesh #1 in Table 1) and 58–59 for the finest mesh (mesh #12 in Table 1), and the relative residuals of the returned iterations were in the range of $(1.1–9.7) \cdot 10^{-10}$.

*4.4. Assessment*

Our main objective was to experimentally investigate superconvergence in the sense of Theorem 1 on the five test problems described in Section 3. To achieve this, we evaluated the order of convergence (in the $H^1$ norm) of the linear part of the discrete velocity to the piecewise-linear nodal interpolation of the exact velocity, and the order of convergence (in the $L^2$ norm) of the discrete pressure to the exact pressure. Since in $H_0^1(\Omega)$ the semi-norm $|*|_{H^1}$ is itself a norm and is equivalent to the usual norm $\|*\|_{H^1}$, either the semi-norm $|*|_{H^1}$ or the norm $\|*\|_{H^1}$ can equivalently be used to assess convergence rates in the sense of Theorem 1. Here, we used the norm $\|*\|_{H^1}$.

In addition, we also assessed how well the piecewise-linear part of the computed velocity $\underline{u}_{hl}$ and the complete computed velocity $\underline{u}_h = \underline{u}_{hl} + \underline{u}_{hb}$ approximate the exact velocity $\underline{u}$. As previously noted, in practical applications $\underline{u}_{hl}$ has often been observed to approximate $\underline{u}$ better than $\underline{u}_h$ thereby suggesting that the bubble function, which is crucial to stabilize the MINI element, has no noticeable effect on the quality of the velocity approximation [12–16]. Most of the empirical observations that back up this result are however restricted to the two-dimensional case, hence the interest for the present assessment in three dimensions. Operatively, besides comparing the respective approximating errors in the H$^1$ and L$^2$ norms (i.e. $\| \underline{u} - \underline{u}_h \|_{H^1}$ versus $\| \underline{u} - \underline{u}_{hl} \|_{H^1}$ and $\| \underline{u} - \underline{u}_h \|_{L^2}$ versus $\| \underline{u} - \underline{u}_{hl} \|_{L^2}$), we also compared the respective L$^2$ norms of the divergence (i.e. $\| div(\underline{u}_h) \|_{L^2}$ versus $\| div(\underline{u}_{hl}) \|_{L^2}$). Conservation of mass requires the exact velocity $\underline{u}$ to be divergence-free, which in the weak formulation reduces to $\| div(\underline{u}) \|_{L^2} = 0$. Low-order mixed finite elements such as the MINI element are only asymptotically divergence-free, meaning that $\| div(\underline{u}_h) \|_{L^2}$ converges to zero as the mesh is gradually refined. Therefore, the comparison between $\| div(\underline{u}_h) \|_{L^2}$ and $\| div(\underline{u}_{hl}) \|_{L^2}$ is informative to assess which one, $\underline{u}_h$ or $\underline{u}_{hl}$, conserves mass better.

## 5. Results and Discussion

Convergence histories and convergence rates are provided in Figure 4 and in Table 2, respectively. As can be noted in Table 2 (first and second rows from the top), the rate of convergence of $\| \underline{u} - \underline{u}_h \|_{L^2}$ is within 2.13–2.19 and that of $\| \underline{u} - \underline{u}_h \|_{H^1}$ is within 1.08–1.18: these rates are consistent with the $O(h^2)$ and $O(h)$ convergence expected from standard finite element theory [9], thereby reassuring on the correctness of the present methodology and of its implementation.

**Table 2.** Experimental convergence rates.

|  | Prob. #1 | Prob. #2 | Prob. #3 | Prob. #4 | Prob. #5 |
|---|---|---|---|---|---|
| $\| \underline{u} - \underline{u}_h \|_{L^2}$ | 2.15 | 2.19 | 2.17 | 2.13 | 2.13 |
| $\| \underline{u} - \underline{u}_h \|_{H^1}$ | 1.08 | 1.11 | 1.09 | 1.18 | 1.11 |
| $\| P - P_h \|_{L^2}$ | 1.48 | 1.61 | 1.49 | 2.04 | 1.65 |
| $\| i_h \underline{u} - \underline{u}_{hl} \|_{L^2}$ | 2.19 | 2.34 | 2.21 | 1.99 | 2.06 |
| $\| i_h \underline{u} - \underline{u}_{hl} \|_{H^1}$ | 1.31 | 1.35 | 1.34 | 1.21 | 1.00 |
| $\| \underline{u} - \underline{u}_{hl} \|_{L^2}$ | 2.15 | 2.19 | 2.17 | 2.13 | 2.12 |
| $\| \underline{u} - \underline{u}_{hl} \|_{H^1}$ | 1.07 | 1.08 | 1.07 | 1.11 | 1.09 |

The rate of convergence of $\| i_h \underline{u} - \underline{u}_{hl} \|_{H^1}$ (fifth row from the top in Table 2) is within 1.00–1.35 which appears consistent with $O(h)$ convergence, indicating that the $O(h^{3/2})$ superconvergence for the velocity in the sense of Theorem 1 is not observed in the present test cases. No sign of superconvergence in velocity in the L$^2$ norm either: the rate of convergence of $\| i_h \underline{u} - \underline{u}_{hl} \|_{L^2}$ is within 1.99–2.34 that seems consistent with $O(h^2)$ convergence, which is the same order of convergence of $\| \underline{u} - \underline{u}_h \|_{L^2}$ and is in line with existing mixed finite element theory [9]. In contrast, the rate of convergence of $\| P - P_h \|_{L^2}$ (third row from the top in Table 2) is within 1.48–1.65 for test problems #1, #2, #3 and #5, which seems consistent with $O(h^{3/2})$ convergence, and equal to 2.04 for test problem #4 which seems consistent with $O(h^2)$ convergence. Therefore, the $O(h^{3/2})$ superconvergence for the pressure in the sense of Theorem 1 is observed in the present test cases, and test problem #4 indicates that quadratic convergence in pressure is indeed possible.

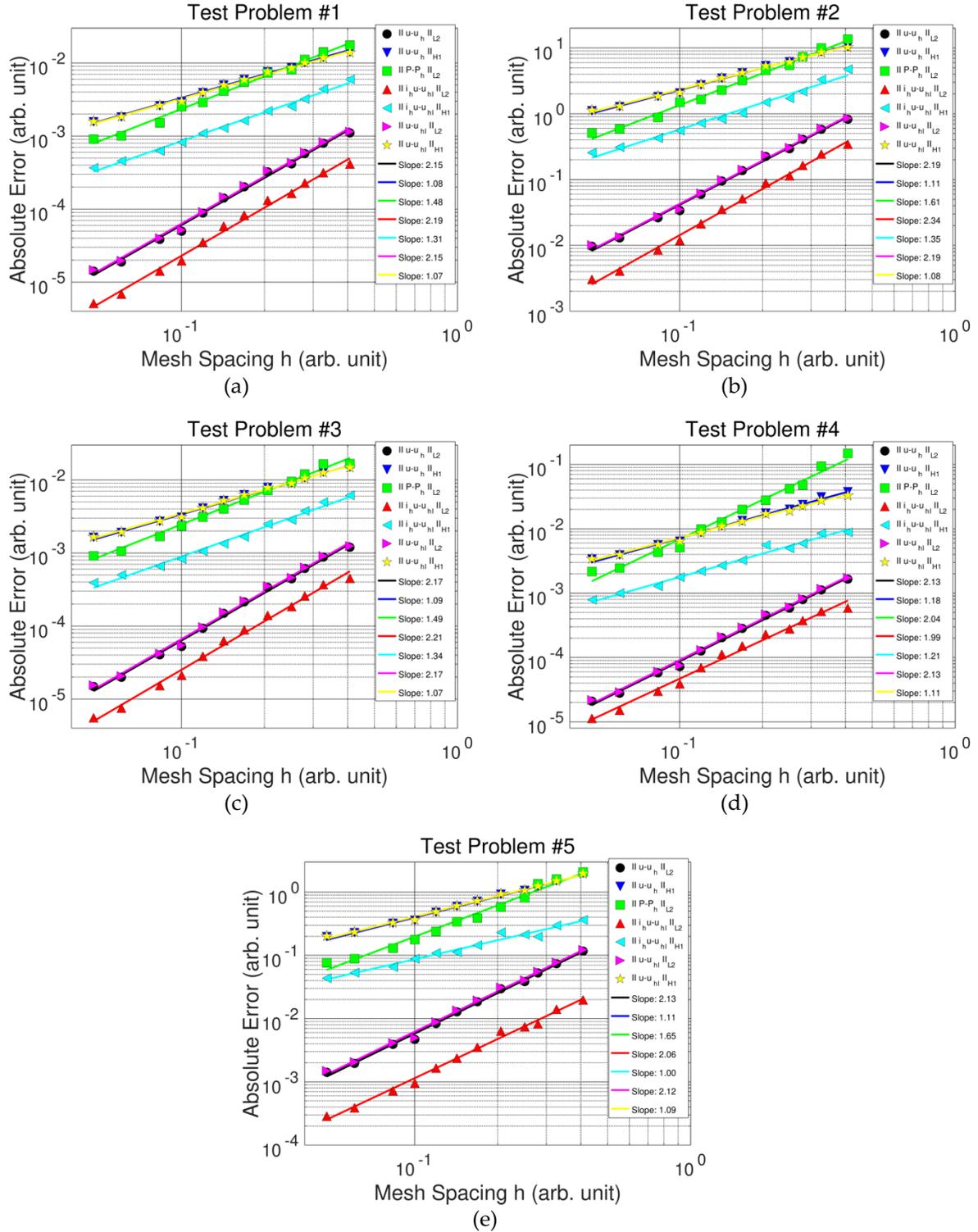

**Figure 4.** Convergence histories for velocity and pressure errors: (a) Test problem #1, (b) Test problem #2, (c) Test problem #3, (d) Test problem #4, (e) Test problem #5 (the solid lines are power law fits through the data points).

Finally, the convergence rate of $\| \underline{u} - \underline{u}_{hl} \|_{L^2}$ is within 2.12–2.19 and that of $\| \underline{u} - \underline{u}_{hl} \|_{H^1}$ is within 1.07–1.11: these rates compare favorably with those of $\| \underline{u} - \underline{u}_h \|_{L^2}$ and $\| \underline{u} - \underline{u}_h \|_{H^1}$ previously discussed and also with the corresponding $O(h^2)$ and $O(h)$ convergence rates expected from standard finite element theory [9].

As previously discussed, the existing $O(h^{3/2})$ superconvergence theory embodied in Theorem 1 is restricted to three-directional triangular meshes in two dimensions. Our previous experimental results [11], which are restricted to the two-dimensional case, indicate that the $O(h^{3/2})$ superconvergence can possibly extend to unstructured triangular meshes. The experimental results documented herein suggest that the $O(h^{3/2})$ superconvergence in pressure may extend to unstructured tetrahedral meshes in three dimensions where, notably, quadratic convergence may also be achievable, whereas the $O(h^{3/2})$ superconvergence in velocity might be restricted to the two-dimensional case. In turn, this indicates that there is ample scope to extend the existing superconvergence theory for the Stokes-MINI problem in two and three dimensions.

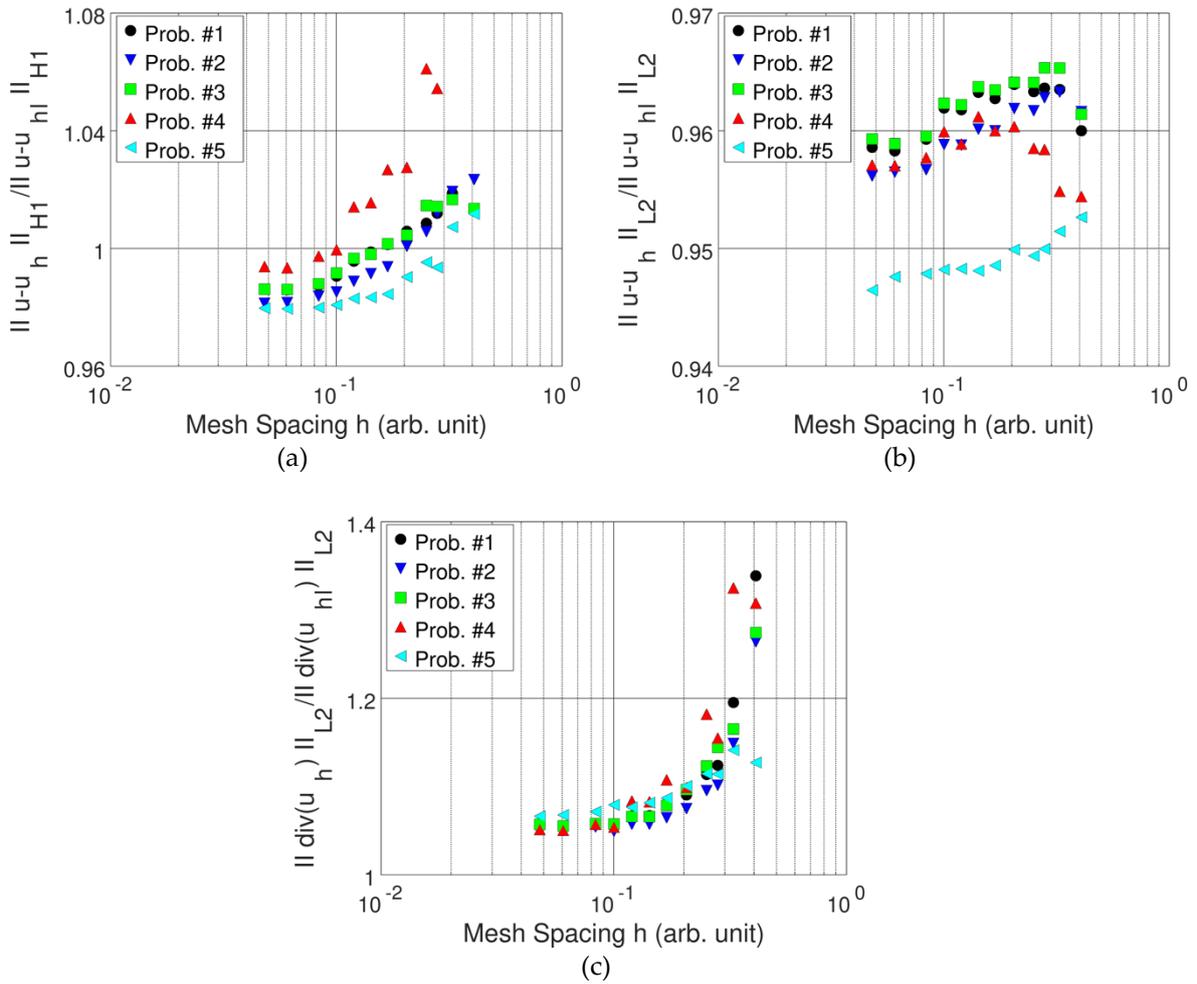

**Figure 5.** Comparison between the complete computed velocity $\underline{u}_h$ and the piecewise-linear part of the computed velocity $\underline{u}_{hl}$: (a) ratio of velocity errors in $H^1$ norm, (b) ratio of velocity errors in $L^2$ norm, (c) ratio of $L^2$ norms of the divergence of the computed velocity.

Regarding the comparison between $\underline{u}_{hl}$ and $\underline{u}_h = \underline{u}_{hl} + \underline{u}_{hb}$, the ratios of the respective errors and of the respective divergence values are provided in Figure 5 as a function of the mesh spacing parameter. As can be noted in Figure 5(a), the ratios of $\|\underline{u} - \underline{u}_h\|_{H^1}$ and $\|\underline{u} - \underline{u}_{hl}\|_{H^1}$ are above one on coarse meshes and decrease below one when the mesh is sufficiently refined, thereby indicating that $\underline{u}_h$

approximates $\underline{u}$ in the H[1] norm better than $\underline{u}_{hl}$ only when the mesh is fine enough whilst the opposite is true on coarse meshes, and the distinction between coarse and fine meshes is problem-dependent. On the other hand, the ratios of $\| \underline{u} - \underline{u}_h \|_{L^2}$ and $\| \underline{u} - \underline{u}_{hl} \|_{L^2}$ in Figure 5(b) are strictly lower than one, indicating that $\underline{u}_h$ approximates $\underline{u}$ in the L² norm better than $\underline{u}_{hl}$ in all tests. Finally, it is evident in Figure 5(c) that $\| div(\underline{u}_h) \|_{L^2}$ is always bigger than $\| div(\underline{u}_{hl}) \|_{L^2}$ thereby indicating that $\underline{u}_{hl}$ conserves mass better than $\underline{u}_h$ in all tests.

In summary, therefore, in all test problems $\underline{u}_h$ approximates the exact velocity in the L² norm better than $\underline{u}_{hl}$, and this is also the case for the H¹ norm provided that the mesh is fine enough. On the other hand, $\underline{u}_{hl}$ conserves mass better than $\underline{u}_h$. As it is evident in Figure 5 these effects are numerically small, particularly so on fine meshes. These results, together with our previous findings in the two-dimensional case [11], suggest that the key role of the bubble function is to stabilize the MINI formulation, whilst its effect on the accuracy of the velocity approximation is minor. This clearly justifies the use of the simpler and easier to compute piecewise-linear part of the computed velocity $\underline{u}_{hl}$ for a posteriori error estimation and in postprocessing.

## 6. Concluding Remarks

Using five purpose-made three-dimensional benchmark test cases with analytical solution, we have experimentally investigated the MINI mixed finite element discretization of the Stokes problem using unstructured tetrahedral meshes. The focus of the investigation was on $O(h^{3/2})$ superconvergence and on assessing how well the complete computed velocity and the piecewise-linear part of the computed velocity approximate the exact velocity. Our main conclusions can be summarized as follows:

- Regarding the pressure error $\| P - P_h \|_{L^2}$, $O(h^{3/2})$ superconvergence was observed in four test problem and $O(h^2)$ superconvergence was observed in one test problem, whereas no superconvergence was observed for the velocity error $\| i_h\underline{u} - \underline{u}_{hl} \|_{H^1}$. Together with the documented observations in the two-dimensional case [11,29,30], these results suggest a far more general validity of the Stokes-MINI superconvergence than what the existing theory covers. This indicates that there is ample scope to extend the existing superconvergence theory for the Stokes-MINI problem in two and three dimensions.

- In all test problems, the exact velocity $\underline{u}$ is well approximated by both the complete computed velocity $\underline{u}_h$ and by its piecewise-linear part $\underline{u}_{hl}$: the former is generally closer to the exact velocity in the L² norm (this is also the case in the H¹ norm when the mesh is fine enough), whereas the latter conserves mass better. These effects are however numerically small, particularly so on fine meshes. Together with our previous findings in the two-dimensional case [11], these results suggest that the key role of the bubble function is to stabilize the MINI formulation, whilst its effect on the quality of the velocity approximation is minor. This clearly justifies the use of the simpler and easier to compute piecewise-linear part $\underline{u}_{hl}$ in postprocessing.

The practical relevance of the purpose-made benchmark test cases with analytical solution documented here goes beyond the present study, particularly when considering the paucity in the open literature of such test cases for the Stokes problem in three dimensions.

**Author Contributions:** Conceptualization, A.C. and D.B.; methodology, A.C. and D.B.; software, A.C.; formal analysis, A.C.; writing—original draft preparation, A.C.; writing—review and editing, A.C. and D.B.

**Funding:** This research received no external funding.